\documentclass[12pt]{article}
\usepackage{mathrsfs}
\usepackage{epic,eepic,epsf,epsfig}
\usepackage{amsfonts,srcltx,mathrsfs}
\usepackage{multirow}
\usepackage{amsmath}
\usepackage{amssymb}
\usepackage{amsbsy}
\usepackage{graphicx}
\usepackage{amsfonts}
\usepackage{color}
\usepackage{setspace}

\newtheorem{lem}{Lemma}[section]%
\newtheorem{theorem}[lem]{Theorem}%

\def\nd{\mathrel{\bigm|\kern-.7em/}}

\def\f{\noindent}

\def\P\GammaL{\hbox{\rm P\GammaL}}

\begin{document}
\title{A tight upper bound of spectral radius in terms of degree deviation}

\footnotetext{E-mails: zhangwq@pku.edu.cn}

\author{Wenqian Zhang\\
{\small School of Mathematics and Statistics, Shandong University of Technology}\\
{\small Zibo, Shandong 255000, P.R. China}}
\date{}
\maketitle

\begin{abstract}
Let $G$ be a graph with $n$ vertices and $m$ edges. The spectral radius $\rho(G)$ of $G$ is the largest eigenvalue of the adjacency matrix of $G$. As is well known, $\rho(G)\geq\frac{2m}{n}$ with equality if and only if $G$ is regular. To bound $\rho(G)-\frac{2m}{n}$,  Nikiforov (2006) introduced the degree deviation of $G$ as $$s(G)=\sum_{1\leq i\leq n}|d_{i}-\frac{2m}{n}|,$$
where $d_{1},d_{2},\ldots,d_{n}$ are the degrees of the vertices of $G$. Nikiforov conjectured that  $\rho(G)-\frac{2m}{n}\leq\sqrt{\frac{1}{2}s(G)}$ for sufficiently large $m$ and $n$.
In this paper, we settle this conjecture without the assumption that $m$ and $n$ are large.
\bigskip

\f {\bf Keywords:} spectral radius; degree deviation; irregularity; adjacency matrix.

\end{abstract}

\baselineskip 17 pt

\section{Introduction}

All graphs considered here are finite, undirected and without loops or multiple edges. For a graph $G$, let $V(G)$ and $E(G)$ denote the vertex set and edge set of $G$, respectively. Let $1,2,...,n$ be all the vertices of $G$. The {\em adjacency matrix} $A(G)=(a_{ij})$ of $G$ is an $n\times n$ square matrix indexed by the vertices of $G$, where $a_{ij}=1$ if and only if $i$ is adjacent to $j$ in $G$. The {\em spectral radius} $\rho(G)$ of $G$ is the largest eigenvalue of $A(G)$. The study of spectral radius of graphs has attracted a large amount of researchers (for example, see \cite{B,CG,HS,LW,N1,R,SW}). For any terminology used but not defined here, one may refer to \cite{BH,CRS}.

Let $G$ be a graph with $n$ vertices and $m$ edges. Collatz and Sinogowitz \cite{CS} showed that $\rho(G)\geq\frac{2m}{n}$. Since the
equality holds if and only if $G$ is a $\frac{2m}{n}$-regular graph,  they proposed the value $\rho(G)-\frac{2m}{n}$  as a
relevant measure of irregularity of $G$.  To bound $\rho(G)-\frac{2m}{n}$, Nikiforov \cite{N} introduced the {\em degree deviation} of $G$ as $$s(G)=\sum_{u\in V(G)}|d_{i}-\frac{2m}{n}|,$$
 where $d_{1},d_{2},\ldots,d_{n}$ are the degrees of the vertices of $G$. For the upper bound of $\rho(G)-\frac{2m}{n}$, Nikiforov \cite{N} proposed the following conjecture.
 
 \medskip
 
 \f{\bf Conjecture 1} \rm{(Nikiforov \cite{N})}. Let $G$ be a graph with $n$ vertices and $m$ edges. Then $\rho(G)-\frac{2m}{n}\leq\sqrt{\frac{1}{2}s(G)}$ for sufficiently large $m$ and $n$.

\medskip

 The constant $\frac{1}{2}$ in Conjecture 1 is best possible. For example, let $G$ be the star graph on $n$ vertices. Note that $m=|E(G)|=n-1$ and $\rho(G)=\sqrt{n-1}$. It is easy to check that
$$\frac{\rho(G)-\frac{2m}{n}}{\sqrt{s(G)}}\rightarrow\sqrt{\frac{1}{2}}$$
when $n\rightarrow\infty$. 

For the progress on  Conjecture 1, (without the assumption that $n$ and $m$ are large), Nikiforov \cite{N}  proved that $\rho(G)-\frac{2m}{n}\leq\sqrt{s(G)}$.  Zhang \cite{Z} proved that $\rho(G)-\frac{2m}{n}\leq\sqrt{\frac{9}{10}s(G)}$. Very recently, Rautenbach and Werner \cite{RW} showed that $\rho(G)-\frac{2m}{n}\leq\sqrt{\frac{2}{3}s(G)}$. In this paper, we settle Conjecture 1 in a stronger form as follows.

\begin{theorem}\label{main1}
Let $G$ be a graph with $n$ vertices and $m$ edges.  Then $\rho(G)-\frac{2m}{n}\leq\sqrt{\frac{1}{2}s(G)}$.
\end{theorem}

To begin the proof of Theorem \ref{main1}, we point out for the readers that  Rautenbach and Werner  obtained their result \cite{RW} before this paper, and the method used there is different from the current one.

\section{The Proof of Theorem \ref{main1} }

For a matrix $A$, let $r_{i}(A)$ be the $i$-th row sum of $A$. The following lemma is deduced from  Lemma 2.1 of \cite{EZ}.

\begin{lem}\rm (Ellingham and  Zha \cite{EZ})\label{polynomial}
Let $G$ be a graph with spectral radius $\rho$, and let $A$ be its adjacency matrix. For any polynomial $f(x)$, we have
$$f(\rho)\leq \max_{u\in V(G)}r_{u}(f(A)).$$
\end{lem}

Now we are in the stage to prove the main theorem of this paper.

\medskip

\f{\bf Proof of Theorem \ref{main1}}.
Let $d=\lceil\frac{2m}{n}\rceil$. For any $u\in V(G)$, let $N(u)$ be the set of vertices adjacent to $u$ in $G$, and let $d(u)=|N(u)|$. For any integer $k\geq0$, let $W_{\geq k}$ (and $W_{\leq k}$, respectively) be the set of vertices of degree at least $k$ (and at most $k$, respectively) in $G$. Let $$s=s(G)=\sum_{v\in V(G)}|d(v)-\frac{2m}{n}|.$$
 Since $\sum_{v\in V(G)}(d(v)-\frac{2m}{n})=0$, we have $$\sum_{v\in W_{\geq d}}(d(v)-\frac{2m}{n})=\frac{1}{2}s.$$
  Thus $$\sum_{v\in W_{\geq d}}(d(v)-d)\leq\sum_{v\in W_{\geq d}}(d(v)-\frac{2m}{n})=\frac{1}{2}s.$$

Let $A=A(G)$ and $\rho=\rho(G)$. We will show that $r_{u}(A^{2}-(d-1)A)\leq d+\frac{1}{2}s$ for any $u\in V(G)$ in two cases as follows.

\f{\bf Case 1}. $d(u)\leq d-1$.

Note that $u\notin W_{\geq d}$ in this case.
Then
 \begin{equation}
\begin{aligned}
r_{u}(A^{2})&=\sum_{v\in N(u)}d(v)\\
&=\left(\sum_{v\in N(u)\cap W_{\leq d-1}}d(v)\right)+\left(\sum_{v\in N(u)\cap W_{\geq d}}d(v)\right)\\
&\leq(d-1)|N(u)\cap W_{\leq d-1}|+\left(\sum_{v\in N(u)\cap W_{\geq d}}(d(v)-d)\right)+d|N(u)\cap W_{\geq d}|\\
&\leq d|N(u)|+\sum_{v\in W_{\geq d}}(d(v)-d)\\
&\leq (d-1)d(u)+d-1+\frac{1}{2}s\\
&<(d-1)d(u)+d+\frac{1}{2}s.
\end{aligned}\notag
\end{equation}
Note that $r_{u}(A)=d(u)$. It follows that $r_{u}(A^{2}-(d-1)A)\leq d+\frac{1}{2}s$.

\f{\bf Case 2}. $d(u)\geq d$.

 Note that $u\in W_{\geq d}$ in this case.
Then
 \begin{equation}
\begin{aligned}
r_{u}(A^{2})&=\sum_{v\in N(u)}d(v)\\
&=\left(\sum_{v\in N(u)\cap W_{\leq d-1}}d(v)\right)+\left(\sum_{v\in N(u)\cap W_{\geq d}}d(v)\right)\\
&\leq(d-1)|N(u)\cap W_{\leq d-1}|+\left(\sum_{v\in N(u)\cap W_{\geq d}}(d(v)-d)\right)+d|N(u)\cap W_{\geq d}|\\
&\leq d|N(u)|+\sum_{v\in W_{\geq d},v\neq u}(d(v)-d)\\
&=(d-1)d(u)+d(u)+\left(\sum_{v\in W_{\geq d}}(d(v)-d)\right)-(d(u)-d)\\
&= (d-1)d(u)+d+\left(\sum_{v\in W_{\geq d}}(d(v)-d)\right)\\
&\leq(d-1)d(u)+d+\frac{1}{2}s.
\end{aligned}\notag
\end{equation}
It follows that $r_{u}(A^{2}-(d-1)A)\leq d+\frac{1}{2}s$.

Let $f(x)=x^{2}-(d-1)x$. By the discussion of Case 1 and Case 2, we have that $r_{u}(f(A))\leq d+\frac{1}{2}s$ for any $u\in V(G)$. By Lemma \ref{polynomial}, we have $f(\rho)\leq d+\frac{1}{2}s$. That is $\rho^{2}-(d-1)\rho\leq d+\frac{1}{2}s$. It follows that $$\rho\leq\frac{d-1}{2}+\sqrt{\frac{(d+1)^{2}}{4}+\frac{1}{2}s}\leq\frac{d-1}{2}+
\left(\sqrt{\frac{(d+1)^{2}}{4}}+\sqrt{\frac{1}{2}s}\right)=d+\sqrt{\frac{1}{2}s}.$$
Then $\rho-\frac{2m}{n}\leq\rho-d+1\leq1+\sqrt{\frac{1}{2}s}$.

As above, we have shown that $\rho(G)-\frac{2m}{n}\leq1+\sqrt{\frac{1}{2}s(G)}$ for any graph $G$ with $n$ vertices and $m$ edges.

For any integer $t\geq1$, as in \cite{N}, let $G^{(t)}$ be the graph obtained from $G$ by replacing each vertex $u$ in $G$ with an independent set $V_{u}$ of $t$ new vertices, and connecting $x\in V_{u}$ to $y\in V_{v}$ if and only if $u$ is adjacent to $v$ in $G$. By Theorem 5 of \cite{N}, we have $\rho(G^{(t)})=t\rho(G)$. Clearly, the average degree of $G^{(t)}$ equals  $\frac{2m}{n}t$, and $s(G^{(t)})=t^{2}s(G)$. Applying the result in the last paragraph to $G^{(t)}$, we have $t\rho(G)-\frac{2m}{n}t\leq 1+\sqrt{\frac{1}{2}t^{2}s(G)}$. By letting $t\rightarrow\infty$, we obtain that $\rho(G)-\frac{2m}{n}\leq\sqrt{\frac{1}{2}s(G)}$. This completes the proof.\hfill$\Box$

\medskip

\f{\bf Data availability statement}

\medskip

There is no associated data.

\medskip

\f{\bf Declaration of Interest Statement}

\medskip

There is no conflict of interest.

\end{document}